\documentclass[12pt]{article}
\usepackage[latin1]{inputenc}
\usepackage{amssymb,amsmath,amsfonts}
\newtheorem{theorem}{Theorem}
\newtheorem{proposition}[theorem]{Proposition}
\newtheorem{lemma}[theorem]{Lemma}

\newtheorem{remark}[theorem]{Remark}
\newtheorem{remarks}[theorem]{Remarks}

\newcommand{\R}{\mathbb{R}}
\newcommand{\Q}{\mathbb{Q}}
\newcommand{\Sf}{\mathbb{S}}

\newcommand{\Hy}{\mathbb{H}}

\newcommand{\po}{{\hspace*{-1ex}}{\bf .  }}

\newcommand{\nab}{\tilde\nabla}

\newcommand{\End}{\mbox{End}}
\newcommand{\Les}{\mathbb{L}}
\def\<{{\langle}}
\def\>{{\rangle}}

\def\n{\nabla}
\def\d{\partial}
\def\e{\epsilon}
\def\be{\begin{equation} }
\def\ee{\end{equation} }
\def\p{\partial}

\def\proof{\noindent\emph{Proof: }}
\def\qed{\ifhmode\unskip\nobreak\fi\ifmmode\ifinner
\else\hskip5 pt\fi\fi\hbox{\hskip5 pt\vrule width4 pt
height6 pt depth1.5 pt\hskip 1pt}}
\begin{document}

\title{Hypersurfaces of space forms carrying a\\ totally geodesic 
foliation
\thanks {\it Mathematics Subject Classification (2000) 53B25, 53C12, 53C40}
}
\author{ Marcos Dajczer and Ruy Tojeiro\footnote{
Partially supported by  FAPESP grant
2016/23746-6.}}
\date{}
\maketitle

\begin{abstract}
In this paper we  give a complete local parametric classification of the 
hypersurfaces  with dimension at least three of a space form  that carry a 
totally geodesic foliation of codimension one. A classification under the 
assumption that the leaves of the foliation are complete was already given 
in \cite{drt} for Euclidean hypersurfaces. We prove that there exists exactly 
one further class of local examples in Euclidean space, all of which have rank 
two. We also extend the classification under the global assumption 
of completeness of the leaves for hypersurfaces of the sphere and show that 
there exist plenty of examples in hyperbolic space. \vspace{1ex}
\end{abstract}

In \cite{drt} the following basic problem was addressed: Which are the 
Euclidean hypersurfaces  of dimension at least three that carry  
a totally geodesic foliation of codimension one?
Recall that a  smooth foliation ${\cal F}$ on a Riemannian manifold $M^n$ 
is \emph{totally geodesic} if all leaves of ${\cal F}$ are totally geodesic 
submanifolds of $M^n$, that is, if any geodesic of $M^n$ that is tangent 
to ${\cal F}$ at some point  is (locally) contained in the leaf of ${\cal F}$ 
through that point. 

The main result in \cite{drt} states that if $f\colon  M^n\to\R^{n+1}$, 
$n\geq 3$, is an isometric immersion  of a nowhere flat  connected 
Riemannian manifold  that carries a totally geodesic foliation of 
codimension one with \emph{complete} leaves, then  $f$ is either ruled or 
a partial tube over a curve, provided that $f(M)$ does not contain any 
surfacelike strip.  

A hypersurface $f\colon M^n\to\R^{n+1}$, $n\geq 2$, is \emph{ruled} if it 
carries a smooth  foliation of codimension one by (open subsets of) affine 
subspaces of $\R^{n+1}$, whereas it is a  \emph{partial tube} over a smooth 
regular curve $\gamma\colon I\subset\R\to\R^{n+1}$ if it is generated by parallel 
transporting along $\gamma$, with respect to the normal connection, a 
hypersurface $N^{n-1}$ of the (affine) normal space $N_\gamma I(t_0)$ to 
$\gamma$ at some point $t_0\in I$ (see Section $1.1$ for further details).   
If $n\geq 3$, by a \emph{surfacelike strip} we mean an open  subset 
$U\subset M^n$ isometric to a product $L^2\times\R^{n-2} $ along which $f$ 
splits as $f=f_1\times id$, with $f_1\colon L^2\to \R^3$  an isometric immersion 
and $id\colon\R^{n-2}\to\R^{n-2}$  the identity map. 

The goal of this paper is to give a full classification of the  
hypersurfaces that carry a totally geodesic foliation of codimension one 
without requiring the leaves to be complete. A rather weak necessary 
condition  was provided by Theorem $5$ in \cite{drt}. Hence, until now 
the local classification problem remained wide open. 

   We prove that, apart from ruled hypersurfaces, surfacelike hypersurfaces
and partial tubes over curves, all other examples without flat points 
have rank two and can be completely described in terms of the Gauss 
parametrization (see Section $1.2$) by means of surfaces $g\colon L^2\to\Sf^n$ 
in the unit sphere that we call of \emph{type $D$}. These are the surfaces 
that carry a nowhere orthogonal real conjugate  system of coordinates 
$(u,v)$ for which one of the Christoffel symbols vanishes everywhere. 
More precisely, all of the coordinate functions  $g^1,\ldots,g^{n+1}$ of 
$g$ in $\R^{n+1}$ are solutions of the differential equation of second 
order
$$
g^j_{uv}+ag^j_v+bg^j=0
$$
for some $a,b\in C^\infty(L)$, satisfying the constraint 
$\sum_{j=1}^{n+1}(g^j)^2=1$.

  Surfaces of type $D$ are closely related to the surfaces of type $C$ 
that were used  in \cite{df} to give a parametrization of all local 
isometric immersions of $\R^n$ into $\R^{n+2}$.  Inspired by a comment 
in \cite{ff} about surfaces of type $C$, we  show that any surface of type 
$D$ in $\Sf^n$ can be generated by a unit principal direction of a  
surface with flat normal bundle in $\R^{n+1}$ that satisfies some generic 
regularity conditions. In particular, this also shows that surfaces of type 
$D$ come in pairs of dual surfaces.

We recall that all Euclidean surfaces with flat normal bundle can be given 
explicitly in terms of a set of solutions of a completely integrable first 
order linear system of PDEs associated to the vectorial Ribaucour transformation 
as shown in \cite{fe} and \cite{dft}. In particular, since these do not satisfy
any constraint,  there exists an abundance of surfaces 
of type $D$.

The other goal of this paper is to study the problem for hypersurfaces of a 
nonflat space form. In this case, we are also able to give a complete solution 
of the local version of the problem,  as well as of the global version in the 
spherical case. As for the global version in hyperbolic space, there  exist 
plenty of examples with complete leaves, which include the isometric immersions 
$f\colon \Hy^n\to \Hy^{n+1}$, the ruled hypersurfaces and those called 
generalized cones (see Section $4$).  To characterize for which of the 
remaining examples all leaves of the foliation are complete is an open problem.

\section[Preliminaries]
{Preliminaries}

In this section we recall several results on submanifold theory in space  
forms that will be used in the sequel.

\subsection[Partial tubes]
{Partial tubes}

We first recall the precise definition of  partial tubes over  curves, and 
then state a result from \cite{to} which implies that these are precisely 
the solutions to our problem for which the totally geodesic distribution 
is invariant by the shape operator of the hypersurface.\vspace{1,5ex}

Let $\gamma\colon I\subset\R\to\R^{n+1}$ be a smooth unit speed curve 
defined on an interval and let $\{\xi_1,\ldots,\xi_n\}$ be an orthonormal 
set of parallel normal vector fields along $\gamma$. Consider the parallel 
vector bundle isometry   $\phi\colon I\times\R^n\to N_\gamma I$ onto the 
normal bundle of $\gamma$ given by
$$
\phi_{s}(y)=\phi(s,y)=\sum_{i=1}^ny_i\xi_i(s)
$$ 
for all $s\in I$ and $y=(y_1,\ldots,y_n)\in\R^n$. 
Let $f_0\colon M^{n-1}_0\to\R^n$  be a hypersurface, denote 
$M^n=M_0^{n-1}\times I$ and define a map $f\colon M^n\to\R^{n+1}$ by
$$
f(x,s)=\gamma(s)+\phi_s(f_0(x)).
$$
One can check that $f$ is an immersion whenever  
$f_0(M_0)\subset\Omega(\gamma;\phi)$, where
$$
\Omega(\gamma;\phi)=\{Y\in\R^n:\<\gamma''(s),\phi_{s}(Y)\>
\neq 1\;\mbox{for any}\;s\in I\}.
$$
In this case, the hypersurface $f$ is called the \emph{partial tube over 
$\gamma$  with fiber $f_0$}. Geometrically, it is generated 
by parallel transporting $f_0(M_0)$ with respect to the normal connection 
of $\gamma$.
\vspace{1ex}

This construction can be extended as follows to the case in which the 
ambient space is either  the sphere or the hyperbolic space. \vspace{1ex}

Let $\R_\mu^{n+2}$ denote either Euclidean 
space $\R^{n+2}$ or Lorentzian space $\Les^{n+2}$, depending on whether 
$\mu=0$ or $1$, respectively. Denote by $\Q_\e^{n+1}\subset \R_\mu^{n+2}$, 
$\e=1-2\mu$,  either the sphere $\Sf^{n+1}$ or the hyperbolic space $\Hy^{n+1}$.
Let $\gamma\colon I\to\R_\mu^{n+2}$ be a smooth unit speed curve such that 
$\gamma(I)$ is contained in $\Q_\e^{n+1}\subset\R_\mu^{n+2}$.  
Set $\gamma=i\circ\tilde\gamma$ where $\tilde\gamma\colon I\to\Q_\e^{n+1}$
and $i\colon\Q_\e^{n+1}\to\R_\mu^{n+2}$ denotes the inclusion.  
Let $\{\xi_1,\ldots,\xi_n\}$ be a parallel orthonormal frame of the normal 
bundle $N_{\tilde \gamma}I$ of $\tilde{\gamma}$. 
Then  $\{\xi_1,\ldots,\xi_n,\xi_{n+1}=\gamma\}$ is a parallel orthonormal 
frame of the normal bundle $N_{\gamma}I$ of ${\gamma}$ and we may define a  
parallel vector bundle isometry 
$\phi\colon I\times\R_\mu^{n+1}\to N_{\gamma}I$ by
$$
\phi_{s}(y)=\phi(s,y)=\sum_{i=1}^{n+1}y_i\xi_i(s).
$$
Let $e\in\R_\mu^{n+1}$ be such that $\gamma(s)=\phi_{s}(e)$ for all 
$s\in I$, and let $f_0\colon M_0^{n-1}\to\R_\mu^{n+2}$ be an isometric 
immersion such that 
$$
f_0(M_0)\subset\Q_\e^{n+1}\cap\left(e+\Omega(\gamma,\phi)\right)
\subset\R_\mu^{n+2}.
$$
Then the map $f\colon M^n=M_0^{n-1}\times I\to\Q_\e^{n+1}\subset\R_\mu^{n+2}$
defined by 
$$
f(x,s)=\phi_{s}(f_0(x))
$$
is called the \emph{partial tube over $\gamma$ with fiber $f_0$}.
\vspace{1ex}

For a partial tube as above, and endowing $M^n=M_0^{n-1}\times I$ with the induced 
metric, the distribution on $M^n$ given by the tangent spaces to the first factor 
is totally geodesic. Equivalently, the induced metric on $M^n$ is a twisted product 
metric
$$
ds^2=ds_0^2+\rho^2 dt^2
$$
for some $\rho\in C^{\infty}(M)$, where $ds^2_0$ and $dt^2$ are the metrics on 
$M_0^{n-1}$ and $I$, respectively.  Moreover, the tangent vector $\d/\d t$ to the 
second factor is a principal direction of $f$ at any point of $M^n$ (see \cite{to} 
for details).  It follows from more general results given  in \cite{to} 
(see Theorem $16$ and Corollary $18$) that this property characterizes 
partial tubes over curves among hypersurfaces of  $\Q_\e^{n+1}$
that carry a smooth totally geodesic distribution of codimension one.

\begin{proposition}\po\label{pt2} Let $f\colon  M^n\to\Q_\e^{n+1}$ 
be an isometric immersion of a twisted product $M^n=M_0^{n-1}\times_{\rho} I$, 
where $I\subset\R$ is an open interval and $\rho\in C^{\infty}(M)$. 
If the tangent vector $\d/\d t$ to the second factor is a principal direction 
of $f$ at any point of $M^n$, then $f$ is  a partial tube over a curve.
\end{proposition}

\subsection[The Gauss parametrization]
{The Gauss parametrization}

In this section,  we recall from \cite{dg} or \cite{dt} the 
\emph{Gauss parametrization} of an oriented 
hypersurface with constant index of relative nullity in a space form. 
\vspace{1ex}

Given an oriented hypersurface $f\colon M^n\to\Q_\e^{n+1}$ in a simply connected
complete space form of sectional curvature $\e=0,1$ or $-1$, we denote by  $\Delta(x)$ 
the kernel of its shape operator $A$ with respect to a unit normal vector field 
$N$ at any $x\in M^n$, called the \emph{relative nullity} subspace. The
\emph{index of relative nullity} at $x\in M^n$ is $\nu(x)=\dim\Delta(x)$.
\vspace{1ex}

We start with the  case of Euclidean hypersurfaces. 

\begin{proposition}\po
Let $g\colon L^{n-k}\to\Sf^n$ be an isometric immersion of a 
Riemannian manifold and let \mbox{$\gamma\in C^\infty(L)$}. 
Set $h=i\circ g$ where $i\colon \Sf^n\to\R^{n+1}$ is the inclusion map,   
and consider the map $\psi\colon\Lambda\to\R^{n+1}$ 
defined on the normal bundle $\Lambda=N_gL$ of $g$ by
\be\label{gpara}
\psi(y,w)=\gamma(y)h(y)+h_*\,\textup{grad}\,\gamma(y)+i_*w.
\ee
Then, on the open subset of regular points,  $\psi$ is an immersed hypersurface 
with constant index of relative nullity $\nu=k$. 

Conversely, any  hypersurface of $\R^{n+1}$ with constant index of 
relative nullity $\nu=k$ can be parametrized in this way, at least locally. 
The parametrization is global if the leaves of the relative nullity distribution 
are complete.
\end{proposition}
 
According to the orthogonal decomposition $TM=\Delta\oplus\Delta^{\perp}$
we decompose $X\in TM$ as $X=X^v+X^h$. Recall that the \emph{splitting tensor} 
$C\colon \Gamma(\Delta)\times \Gamma(\Delta^{\perp})\to \Gamma(\Delta^{\perp})$ 
is defined as
$$
C(T,X)=-(\nabla_XT)^h
$$
for all $T\in \Gamma(\Delta)$ and $X\in \Gamma(\Delta^{\perp})$. If $x\in M^n$ and 
$T\in\Delta(x)$, then the tensor gives rise to a map
$C_T\colon\Delta^\perp(x)\to\Delta^\perp(x)$. Then we regard $C$ as a map 
$C\colon\Gamma(\Delta)\to\Gamma(\End(\Delta^\perp))$.
\vspace{1ex}

We now collect several properties of the Gauss parametrization needed 
in the sequel.

\begin{proposition}\po\label{gaussprop} 
The following assertions hold:
\begin{itemize}
\item[(i)] The map $\psi$ is regular at $(y, w)\in \Lambda$ if and only if 
the selfadjoint operator  
\be\label{eq:pw}
P_w(y)=\gamma(y)I+\textup{Hess}\,\gamma(y)-A_w
\ee 
on $T_yL$ is nonsingular, where $A_w$ is the shape operator of $g$ with 
respect to $w$. 
\item[(ii)] On the open subset $V$ of regular points, $\psi$ is an immersed 
hypersurface having the map $G\colon\Lambda\to\Sf^n$, given by 
$$
G(y,w)=g(y),
$$
as a  Gauss map of  rank $n-k$.
\item[(iii)] For any  $(y,w)\in V$  there exists a map  
$j=j(y,w)\colon T_yL\to T_{(y,w)}\Lambda$, which is an isometry  
onto the orthogonal complement $\Delta^\perp(y,w)$ of the relative 
nullity subspace $\Delta(y,w)=N_gL(y)$ of $\psi$ at $(y, w)$, such that 
\be\label{eq:nxijx}
\nabla^V_\xi jX=0
\ee
for all $\xi\in \Delta(y,w)$ and $X\in\mathfrak{X}(L)$, where $\nabla^V$ 
denotes the Levi-Civita connection of the metric on $V$ induced by $\psi$,
and such that the shape operator $A$ of $\psi$ at $(y,w)$ with respect to 
$G$  satisfies
\be\label{eq:sffpsi}
A j=-jP_w^{-1}.
\ee
\item[(iv)] For any  $(y,w)\in V$ the splitting tensor 
$C_\xi\colon\Delta^\perp(y,w)\to\Delta^\perp(y,w)$ of $\Delta$
with respect to $\xi\in \Delta(y,w)=N_gL(y)$ is related to the shape operator 
$A_\xi$ of $g$ at $y$ by
\be\label{cxiaxi}
C_\xi j = jA_\xi P_w^{-1}.
\ee
\item[(v)] The Levi-Civita connections $\nabla^L$ and $\nabla^V$ of the metrics 
$\<\,,\,\>_L$ and $\<\,,\,\>_V$ 
on $L^{n-k}$ and $V$ induced by $g$ and $\psi$, respectively, are related by
\be\label{eq:conngpsi}
\<\nabla^L_{P_w^{-1}X} Y, Z\>_L=\<\nabla^V_{jX}jY, jZ\>_V
\ee
for all $X,Y,Z\in\mathfrak{X}(L)$.
\end{itemize}
\end{proposition}

The Gauss parametrization of oriented hypersurfaces $f\colon M^n\to\Q_\e^{n+1}$ 
with $\nu=k$ reads as follows for $\e\in \{-1,1\}$. 

\begin{proposition}\po\label{gauss3} 
Let $g\colon L^{n-k}\to\Sf_\mu^{n+1}$ be an isometric immersion of 
a Riemannian manifold, where
$$
\Sf^n_\mu=\{x\in\R_\mu^{n+1}\colon\<x,x\>=1\}.
$$
Denote
\be\label{gammae}
\Lambda_\e=\{(y,w)\in\Lambda=N_gL:\<w,w\>=\e\},\;\e=1-2\mu,
\ee
and define $\psi\colon\Lambda_\e\to\Q_\e^{n+1}$   by
\be\label{gpara2}
\psi(y,w)= w.
\ee
Then the following assertions hold:
\begin{itemize}
\item[(i)] On the open subset $V$ of regular points, $\psi$ is an immersed 
hypersurface with constant index of relative nullity $\nu=k$, having the 
map $G\colon\Lambda_\e\to\Sf_\mu^{n+1}$, given by 
$$
G(y,w)=g(y),
$$
as a  Gauss map of rank $n-k$.
\item[(ii)] Conversely, any  hypersurface of $\Q_\e^{n+1}$ with constant index of 
relative nullity $\nu=k$ can be parametrized in this way, at least locally. 
The parametrization is global if the leaves of the relative nullity distribution 
are complete.
\item[(iii)] The map $\psi$ is regular at $(y, w)\in \Lambda_\e$ if and only if the 
shape operator $A_w$ of $g$ at $y\in L^{n-k}$ is nonsingular. 
\item[(iv)] At any  $(y,w)\in V$, there exists an isometry 
$j=j(y,w)\colon T_yL\to\Delta^\perp(y,w)$ onto the orthogonal complement 
of the relative nullity subspace 
$$\Delta(y,w)=\{w'\in N_gL(y)\;:\; \<w', w\>=0\}$$ of $\psi$ at $(y, w)$, such 
that (\ref{eq:nxijx}) holds for all $\xi \in \Delta(y,w)$ and such that
the shape operator $A$ of $\psi$ at $(y,w)$ with respect to $G$ satisfies
$$
A j=jA_w^{-1}.
$$
\item[(v)] For any  $(y,w)\in V$ the splitting tensor 
$C_\xi\colon\Delta^\perp(y,w)\to\Delta^\perp(y,w)$ of $\Delta$
with respect to $\xi\in \Delta(y,w)$ is related to the shape operator 
$A_\xi$ of $g$ at $y$ by
\be\label{cxiaxi2}
C_\xi j = jA_\xi A_w^{-1}.
\ee
\item[(vi)] The Levi-Civita connections $\nabla^L$ and $\nabla^V$ of the metrics 
$\<\,,\,\>_L$ and $\<\,,\,\>_V$ on $L^{n-k}$ and $V$ induced by $g$ and $\psi$, 
respectively, are related by
\be\label{eq:conngpsi2}
\<\nabla^L_{A_w^{-1}X} Y, Z\>_L=\<\nabla^V_{jX}jY, jZ\>_V
\ee
for all $X,Y,Z\in\mathfrak{X}(L)$.
\end{itemize} 
\end{proposition}

\section[Surfaces of type $D$] 
{Surfaces of type $D$}

This section is devoted to introduce a class of surfaces in the unit  
sphere $\Sf^n$ that we call surfaces of type $D$. We will see in the next 
section that these surfaces are precisely  the Gauss images of Euclidean 
hypersurfaces of rank two that carry a totally geodesic foliation with 
codimension one. 
\vspace{1ex}

We say that a surface $g\colon L^2\to\Sf^n$, $n\geq 3$, is a 
\emph{surface of type $D$} if there exist linearly independent  
vector fields  $X,Y\in\mathfrak{X}(L)$ with $\|Y\|=1$ such that
\be\label{eq:typeDcond}
\begin{cases}
(i)\;\;\<X,Y\>\neq 0\\
(ii)\;\;\nabla^g_XY=0\\
(iii)\;\;\alpha_g(X,Y)=0
\end{cases}
\ee
where $\nabla^g$ denotes the Levi-Civita connection of the metric on $L^2$ 
induced by $g$. Notice that conditions $(ii)$ and $(iii)$ together are equivalent
to $\bar{\nabla}_X g_*Y=0$, where $\bar{\nabla}$ is the connection of $\Sf^n$, 
that is, they mean that $Y$ is parallel along $X$ in  $\Sf^n$.\vspace{1ex}

If $(u,v)$ are coordinates on $L^2$  such that  $X$ and $Y$ are collinear
with $\p u$ and $\p v$, respectively, then the preceding conditions can be 
written as  
$$
\begin{cases}
(i')\;\; \<\p u,\p v\>\neq 0\\
(ii')\;\; \nabla^g_{\p u}\p v\in\textup{span}\{\p v\}\\
(iii')\;\; \alpha_g(\p u,\p v)=0.
\end{cases}
$$
The last two conditions are equivalent to the position vector of $g$ in 
$\R^{n+1}$ satisfying the differential equation
\be\label{eq:pos}
g_{uv}+ag_v+bg=0
\ee
for smooth functions $a=a(u,v)$ and $b=b(u,v)=\<g_u,g_v\>$, where $g_u$ 
and $g_v$ denote the partial derivatives of $g$.

Surfaces of type $D$ have in common conditions $(ii)$ and $(iii)$ with the
surfaces of type~$C$ defined in \cite{df} (see also \cite{cd}), for which, 
instead of condition $(i)$, one  requires that the vector field $Y$ be 
nowhere an asymptotic direction for $g$.  
Surfaces of type $C$ were used in \cite{df} to give  a parametrization of 
all local isometric immersions of $\R^n$ into $\R^{n+2}$. 
\vspace{1ex}

The following result shows that the existence of linearly independent  
vector fields  $X,Y\in\mathfrak{X}(L)$ with $\|Y\|=1$  satisfying
conditions $(i)$ to $(iii)$ imposes no restrictions on a surface 
$g\colon L^2\to\Sf^3$ with Gauss curvature $K\neq 1$ everywhere.

\begin{proposition}\po \label{n3} Any surface $g\colon L^2\to\Sf^3$
with Gauss curvature $K\neq 1$ everywhere is a surface of type $D$.
\end{proposition}

\proof Since $K\neq 1$ everywhere, a unit normal vector field to $g$
defines an immersion $\hat g\colon L^2\to\Sf^3$, called the \emph{polar
surface} of $g$. The metrics $\<\,,\,\>_g$ and $\<\,,\,\>_{\hat g}$ 
on $L^2$ induced by $g$ and $\hat g$, respectively, are related by
$$
\<X,Y\>_{\hat g}=\<AX,AY\>_g
$$
for all $X,Y\in\mathfrak{X}(L)$, where $A$ is the shape operator of $g$
with respect to $\hat g$, and the corresponding Levi-Civita connections 
$\nabla^g$ and $\nabla^{\hat g}$ are related by
$$
A\nabla^{\hat g}_XY=\nabla^g_XAY
$$
for all $X,Y\in\mathfrak{X}(L)$.

Now let $X,\tilde Y$ be an orthonormal frame of $L^2$ with 
respect to the metric induced by $\hat g$ such that the integral curves 
of $X$ are geodesics with respect to that metric and such that $X$ is 
nowhere a principal direction of $\hat g$. Thus
$$
\nabla^{\hat g}_X\tilde Y=0\;\;\;
\mbox{and}\;\;\;\<\hat AX,\tilde Y\>_{\hat g}\neq 0
$$
everywhere, where $\hat A=A^{-1}$ is the shape operator of 
$\hat g$ with respect to $g$. If $Y=A\tilde Y$, then
$$
\<Y,Y\>_g=\<A^{-1}Y,A^{-1}Y\>_{\hat g}=\<\tilde Y,\tilde Y\>_{\hat g}=1,
$$
$$
\nabla^g_XY= A\nabla^{\hat g}_XA^{-1}Y=A\nabla^{\hat g}_X\tilde Y=0,
$$
$$
\<AX,Y\>_g=\<X,A^{-1}Y\>_{\hat g}=\<X,\tilde Y\>_{\hat g}=0
$$
and 
$$
\<X,Y\>_g=\<A^{-1}X,A^{-1}Y\>_{\hat g}=\<A^{-1}X, \tilde Y\>_{\hat g}
=\<\hat AX,\tilde Y\>_{\hat g}\neq 0
$$
everywhere.
Thus $X,Y$ satisfy conditions $(i)$ to $(iii)$ for $g$. \vspace{1ex}\qed

The next result shows that, generically, surfaces of type $D$ have a dual 
surface of the same type.  

\begin{proposition}\po\label{dual} Let $g\colon L^2\to\Sf^n$ be a 
surface of type $D$ with Gauss curvature $K\neq 1$ everywhere, and let
$X_g,Y_g$ be vector fields satisfying conditions (\ref{eq:typeDcond}) for $g$.
Then the map $k\colon L^2\to\Sf^n$ given by
$
k=g_*Y_g
$
also defines a surface of type $D$.
\end{proposition}

\proof  We have
\be\label{eq:k1}
k_*X_g=-\<X_g,Y_g\>_g\, g
\ee
and
$$
k_*Y_g=g_*\nabla^g_{Y_g} Y_g+\alpha_g(Y_g, Y_g)-g
$$
by the Gauss formula of $g$. Here $\<\,,\,\>_g$ stands for the 
metric induced by $g$ and $\nabla^g$ by its induced connection.
Since $K\neq 1$ everywhere and 
$\alpha_g(X_g,Y_g)=0$, then $\alpha_g(Y_g,Y_g)\neq 0$ everywhere
by the Gauss equation of $g$. Thus $k$ is an immersion.

Denote $X_k=Y_g$ and $Y_k=X_g/\<X_g, Y_g\>$. Then (\ref{eq:k1}) can 
be written as
\be\label{eq:dualk} 
k_*Y_k=-g.
\ee
In particular, $Y_k$ has unit length with respect to the metric 
$\<\,,\,\>_k$ induced by $k$, and
$$
\<X_k,Y_k\>_k=\<k_*Y_g, -g\>=1.
$$
Moreover,
$$
\tilde{\nabla}_{X_k}k_*Y_k=-\tilde{\nabla}_{X_k}g=-g_*Y_g=-k
$$
where $\tilde{\nabla}$ stands for the derivative in $\R^{n+1}$. 
Thus $\bar{\nabla}_{X_k}k_*Y_k=0$, and hence $X_k,Y_k$ satisfy 
conditions (\ref{eq:typeDcond}) for $k$.\vspace{1ex}\qed

Notice that, by (\ref{eq:dualk}), the dual surface of $k$ is $-g$.
Notice also that the position vector fields in Euclidean space of 
the dual  surfaces are orthogonal. \vspace{1ex}

The next result, inspired by a comment in \cite{ff} about surfaces
of type $C$, shows that any  surface  of type $D$ is spanned by a unit 
principal direction of a surface with flat normal bundle in $\R^{n+1}$ 
that satisfies some generic regularity conditions. 

\begin{theorem}\po\label{typeD} 
Let $h\colon L^2\to\R^{n+1}$, $n\geq 3$, be a surface with flat normal 
bundle and let $X,Y\in\mathfrak{X}(L)$ be an orthonormal frame that 
diagonalizes the second fundamental form of $h$. Assume that 
the integral curves of $X$ and $Y$ are nowhere geodesics and that the 
integral curves of $X$ are nowhere asymptotic.  
Then the map $g\colon L^2\to\Sf^n\subset\R^{n+1}$ given~by 
$$
g=h_*X
$$ 
defines a surface of type $D$. Conversely, any surface of type $D$ can 
be locally constructed in this way.
\end{theorem}

\proof Taking derivatives in $\R^{n+1}$ yields
$$
g_*Y=\nab_Yh_*X=-\<\n^h_Y Y, X\>_hh_*Y,
$$
since $\alpha_h(X,Y)=0$ by assumption. Here $\<\,,\,\>_h$ denotes the 
metric induced by $h$ and $\n^h$ stands for its Levi-Civita connection. 
Moreover, since $\n^h_YY\neq 0$ at any point of $L^2$ by assumption, 
then also $g_*Y\neq 0$ at any point of $L^2$. We have 
$$
g_*X=\nab_Xh_*X=\<\n^h_XX, Y\>_hh_*Y+ \alpha_h(X,X).
$$
Thus the map $g$ is an immersion since, by assumption, $\alpha_h(X,X)\neq 0$ 
at any point of $L^2$. 
Moreover, the vector field 
$$
\hat Y=\frac{1}{\<\n^h_YY,X\>_h}Y
$$
has unit length with respect to the metric induced by $g$, and
$$
\nab_Xg_*\hat Y=-\nab_Xh_*Y=\<\n^h_XX,Y\>_hh_*X=\<\n^h_XX,Y\>_h g. 
$$
Thus conditions $(ii)$ and $(iii)$ are satisfied for the vector fields 
$\hat X=X$ and $\hat Y$. Moreover, since also $\n^h_XX\neq 0$ at any 
point of $L^2$, then
$$
\<g_*\hat X,g_*\hat Y\>\neq 0,
$$
and hence $g$ is of type $D$. \vspace{1ex}

Conversely, let $g\colon L^2\to\Sf^n\subset\R^{n+1}$, $n\geq 3$, be a 
surface of type $D$ and let $(u,v)$ be coordinates in a neighborhood $U$ of
a point of $L^2$ such that conditions $(i')$ to $(iii')$ are satisfied.
Thus, the position vector $g$ satisfies \eqref{eq:pos} on $U$.  
Shrinking $U$, if necessary, one can find, using for instance the Riemann 
method, a nowhere vanishing  solution $\varphi\in C^\infty(U)$ of the 
hyperbolic PDE
\be\label{eq:varphi}
\varphi_{uv}-(a+b_u/b)\varphi_v+b\varphi=0
\ee
such that $\varphi_v$ is also nowhere vanishing on $U$. By assumption, 
$b=\<g_u,g_v\>\neq 0$. Observe that \eqref{eq:varphi} is an equation of the 
same type as \eqref{eq:pos}. Now consider the system of PDEs 
\be\label{eq:system}
\begin{cases}h_u=\varphi g\\
h_v=-b^{-1}\varphi_vg_v
\end{cases}
\ee
for $h\colon U\to\R^{n+1}$. It is easily checked that the integrability 
conditions of \eqref{eq:system} are satisfied by virtue of \eqref{eq:pos} 
and \eqref{eq:varphi}. Since $\<h_u,h_v\>=0$ and 
$$
h_{uv}=\frac{\varphi_v}{\varphi}h_u-\frac{b\varphi}{\varphi_v}h_v\in h_*TU
=\textup{span}\{g,g_v\},
$$
then $h$ is a surface with flat normal bundle and $g=h_*X$, where $X$ is 
the unit vector field in the direction of $\p u$ with respect to the metric 
induced by $h$. 

It remains to check that the coordinate curves are nowhere geodesics with 
respect to the metric on $U$ induced by $h$, and that the $u$-coordinate 
curves are nowhere asymptotic for $h$. In other words, we need to verify 
that at any point of $U$ the following conditions hold.
$$ 
\begin{cases}
(i")\;\;\nabla^h_{\p u}\p u\not\in\textup{span}\{\p u\}\\
(ii")\;\;\nabla^h_{\p v}\p v\not\in\textup{span}\{\p v\}\\
(iii")\;\;\alpha_h(\p u,\p u)\neq 0.
\end{cases}
$$
For $(i")$: we have $\nabla^h_{\p u}\p u\in\textup{span}\{\p u\}$
if and only if the component of 
$$
h_{uu}=\varphi_ug+\varphi g_u
$$ 
in $h_*TU=\textup{span}\{g,g_v\}$ is collinear with $h_u=\varphi g$, which cannot 
happen because $\varphi\neq 0$ and $\<g_u,g_v\>\neq 0$ at any point of $U$.
\vspace{1ex}\\
For $(ii")$: if $\nabla^h_{\p v}\p v\in\textup{span}\{\p v\}$, then the $h_*TU$ 
component of 
$$
h_{vv}=-(b^{-1}\varphi_v)_vg_v-b^{-1}\varphi_vg_{vv}
$$ 
is collinear with $h_v=-b^{-1}\varphi_v g_v$. 
Hence, since $\varphi_v\neq0$ then $0=\<g_{vv},g\>=-\<g_v,g_v\>$, and this 
is a contradiction. \vspace{1ex}\\
For $(iii")$:  if $\alpha_h(\p u,\p u)= 0$ then 
$h_{uu}=\varphi_u g+\varphi g_u\in\textup{span}\{g,g_v\}$, 
and this is again a contradiction. \qed

\begin{remarks}\po\label{remark} {\em
\noindent $(1)$  The definition of surfaces of type $D$ can 
be extended for surfaces in the Lorentzian sphere $\Sf_1^{n+2}$ of 
constant sectional curvature $1$ in the Lorentz space $\Les^{n+3}$.
We will see in the last section that these are the Gauss images of 
hypersurfaces with rank two of $\Hy^{n+1}$ that carry a totally geodesic 
foliation of codimension one.
\vspace{1ex}

\noindent $(2)$ In view of Proposition \ref{n3}, it is natural to 
exclude from the definition of surfaces $f\colon L^2\to \Sf^n$
of type $D$ those with Gauss curvature $K\neq 1$ everywhere 
that are contained in an umbilical sphere
$\Sf^3(c)\subset \Sf^n$ of constant curvature $c\geq 1$. Since it is 
easily checked that Proposition \ref{n3} still holds for surfaces 
in the Lorentzian sphere $\Sf_1^{3}$ with Gauss curvature $K\neq 1$
everywhere, we also exclude from the 
definition of surfaces $f\colon L^2\to \Sf_1^{n+2}$
of type $D$ those with Gauss curvature $K\neq 1$ that are contained 
in an umbilical Lorentzian sphere  $\Sf_1^3(c)\subset \Sf_1^{n+2}$
of constant curvature $c\geq 1$.
\vspace{1ex}

\noindent $(3)$  In the notations of Theorem \ref{typeD} we have 
$$
\nab_Yg_*\hat Y=\nab_Yh_*Y=\<\n^h_YY, X\>_h g+\alpha_h(Y,Y). 
$$
Thus $\hat Y$ is an asymptotic direction of $g$ if and only the vectors  
$\alpha_h(Y,Y)$ and $\alpha_h(X,X)$ are linearly dependent. Therefore $g$ is 
also a surface of type $C$ if the first normal spaces of $h$, that is, the 
subspaces spanned by its second fundamental form, have dimension two at any 
point.
\vspace{1ex}

\noindent $(4)$ Notice that the surface given by Proposition \ref{dual} 
is just the other unit principal direction of $h$ in Theorem \ref{typeD}. 
}\end{remarks}

\section[The case of Euclidean hypersurfaces]
{The case of Euclidean hypersurfaces}

In this section we give the local classification of hypersurfaces 
$f\colon M^n\to\R^{n+1}$ that carry  a  totally geodesic distribution of 
codimension one. 
\vspace{1ex}

First we observe that, besides the ruled hypersurfaces referred to in the 
introduction and the partial tubes described in Section $1.1$, there exist 
three trivial families of examples. 
\vspace{1ex}

\noindent $(1)$ The first family is that of isometric immersions
$f\colon U\to\R^{n+1}$ of  open subsets $U\subset\R^n$, since one may 
consider any  foliation of $U$ by (open subsets of) affine hyperplanes.
\vspace{1ex}

\noindent $(2)$ The second family consists of cylindrical surfacelike 
hypersurfaces. For a surface $g\colon L^2\to\R^3$, let $\mathcal{D}_0$ 
be the one-dimensional distribution on $L^2$ spanned by the tangent 
directions to a foliation of $L^2$ by geodesics.
Set $M^n=L^2\times\R^{n-2}$ and define an isometric immersion
$f\colon M^n\to\R^{n+1}$ by $f=g\times id$, where 
$id\colon\R^{n-2}\to\R^{n-2}$ is the identity map. Then 
$\mathcal{D}=\mathcal{D}_0\oplus \R^{n-2}$ is clearly 
a totally geodesic distribution on $M^n$ of codimension one.
We call $f$ a \emph{cylindrical surfacelike} hypersurface.
\vspace{1ex}

\noindent $(3)$ Similarly, for a surface $g\colon L^2\to\Sf^3\subset\R^4$ 
in the sphere, let $\mathcal{D}_0$ be as before. 
Define  $M^n=L^2\times\R_+\times\R^{n-3}$ 
and $f=C(g)\times id$, where $C(g)\colon L^2\times\R_+\to\R^4$, given by 
$C(g)(x, t)=tg(x)$, is the cone over $g$ in $\R^4$.
Then the distribution $\mathcal{D}=\mathcal{D}_0\oplus \R\oplus \R^{n-3}$ 
on $M^n$ is again totally geodesic of codimension one.
We call $f$ a \emph{conical surfacelike} hypersurface.
\vspace{1ex}

Given a surface $g\colon L^2\to\Sf^n$, $n\geq 3$, and  
$\gamma\in C^{\infty}(L)$, we say that $(g,\gamma)$ is a 
\emph{pair of type $D$} if $g$ is a surface of type $D$ and $\gamma$ 
satisfies the same differential equation (\ref{eq:pos}) as the position 
vector of $g$. 
It follows from Proposition \ref{n3} that if $g\colon L^2\to\Sf^n$
is a surface with Gauss curvature $K\neq 1$ everywhere 
contained in an umbilical  $\Sf^3(c)\subset \Sf^n$, 
$c\geq 1$, and $\gamma=0$, then $(g,\gamma)$ is always a pair of type $D$. 
We show next that this is also the case if $g\colon L^2\to\Sf^n$ is totally
geodesic and $\gamma\in C^{\infty}(L)$ is any smooth function such that
the endomorphism $P=P_\gamma=\gamma I+\textup{Hess}\,\gamma$ is invertible 
everywhere.

\begin{proposition}\po \label{tpairs} If $g\colon L^2\to\Sf^n$ is totally 
geodesic and $\gamma\in C^{\infty}(L)$ is  such that
$P=P_\gamma=\gamma I+\textup{Hess}\,\gamma$ is invertible everywhere,
then $(g,\gamma)$ is a pair of type $D$.
\end{proposition}

\proof Regard $g$ as a map $g\colon L^2\to\Sf^2$ and define 
$f\colon L^2\to\R^3$ by 
$$
f(y)=\gamma(y)h(y) +h_*\textup{grad}\,\gamma.
$$
Then $f_*=h_*P$, the metrics $\<\,,\,\>_h$ and $\<\,,\,\>_f$ on $L^2$ 
induced by $h$ and $f$, respectively, are related by
$$
\<X,Y\>_f=\<PX,PY\>_h
$$
for all $X,Y\in\mathfrak{X}(L)$, the corresponding Levi-Civita connections 
$\nabla^h$ and $\nabla^{f}$ are related by
$$
P\nabla^{f}_X Y=\nabla^h_XPY
$$
for all $X,Y\in\mathfrak{X}(L)$, and the shape operator of $f$ with 
respect to $h$ is
$$
A=P^{-1}.
$$
Let $X, \tilde Y$ be an orthonormal frame of $L^2$ with respect to 
the metric induced by $f$ such that the integral curves of $X$ are 
geodesics with respect to that metric and such that $X$ is nowhere a 
principal direction of $f$.
Thus
$$
\nabla^f_X\tilde Y=0\;\;\;\mbox{and}\;\;\;\<AX,\tilde Y\>_f\neq 0
$$
everywhere. Define $Y=P\tilde Y$. Then
$$
\<Y, Y\>_h=\<P^{-1}Y,P^{-1}Y\>_f= \<\tilde Y,\tilde Y\>_f=1,
$$
$$
\nabla^h_XY=P\nabla^f_XP^{-1} Y=P\nabla^f_X\tilde Y=0,
$$
\be\label{eq:p}
\<PX,Y\>_h=\<X,P^{-1}Y\>_f= \<X,\tilde Y\>_f=0
\ee
and 
$$
\<X, Y\>_h=\<P^{-1}X, P^{-1}Y\>_f= \<P^{-1}X, \tilde Y\>_f
=\<AX, \tilde Y\>_f\neq 0
$$
everywhere.
Thus $X, Y$ satisfy conditions $(i)$ to $(iii)$ for $g$,
and (\ref{eq:p}) implies that $\gamma$ satisfies the
same differential equation (\ref{eq:pos}) as the position vector of $g$.
 \vspace{1ex}\qed

It is therefore natural to exclude from the definition of pairs $(g, \gamma)$ of
type $D$ the cases in which either $g\colon L^2\to\Sf^n$ is a surface 
with Gauss curvature $K\neq 1$ everywhere contained
in an umbilical $\Sf^3(c)\subset \Sf^n$, $c\geq 1$, and $\gamma=0$, or 
$g\colon L^2\to\Sf^n$ is totally geodesic and $\gamma\in C^{\infty}(L)$ 
is any smooth function such that the endomorphism 
$P=P_\gamma=\gamma I+\textup{Hess}\,\gamma$ is invertible 
everywhere. Notice that these are precisely the pairs that give rise,
by means of the Gauss parametrization, to conical and cylindrical
surfacelike hypersurfaces, respectively.\vspace{1ex}

We are now in a position to state the local classification of Euclidean hypersurfaces 
of dimension $n\geq 3$ without flat points that carry a totally geodesic distribution 
of codimension one.

\iffalse
The local classification of Euclidean hypersurfaces of dimension $n\geq 3$ 
without flat points that carry a totally  geodesic distribution of codimension 
one is as follows.
In the next statement, that $(g,\gamma)$ is a \emph{pair of type $D$} means 
that  $g\colon L^2\to\Sf^n$, $n\geq 3$, is a surface of type $D$ and that 
$\gamma\in C^{\infty}(L)$ satisfies the
same differential equation (\ref{eq:pos}) as the position vector of $g$. 
\fi

\begin{theorem}\po\label{thm:main}
Let $f\colon  M^n\to\R^{n+1}$, $n\geq 3$, be an isometric immersion of 
a Riemannian manifold without flat points that carries a totally geodesic  
foliation of codimension one. Then there exists an open and dense subset of 
$M^n$ where  $f$ is locally either (cylindrical or conical) surfacelike, ruled, 
a partial tube over a curve, or  given by the Gauss parametrization by a pair 
$(g,\gamma)$ of type~$D$. Conversely, any of these hypersurfaces carries a 
totally geodesic foliation of codimension one.
\end{theorem}

\iffalse
We observe that  conical surfacelike hypersurfaces are covered by the last 
class of hypersurfaces in the above statement;  see part $(1)$ of Remarks \ref{remark}.
\vspace{1ex}
\fi

The proof of  Theorem \ref{thm:main} relies on the following lemma.

\begin{lemma}\po\label{totgeofol}
Let $g\colon L^{2}\to\Sf^n$, $n\geq 3$, be a surface, 
let \mbox{$\gamma\in C^\infty(L)$} and let
$\psi\colon\Lambda\to\R^{n+1}$ be the map
defined on the normal bundle $\Lambda=N_gL$ of $g$ by
$$
\psi(y,w)=\gamma(y)h(y)+h_*\,\textup{grad}\,\gamma(y)+i_*w,
$$
where  $h=i\circ g$ is the composition
of $g$ with the inclusion  $i\colon \Sf^n\to\R^{n+1}$.  
Assume that  there exist $Y,Z\in\mathfrak{X}(L)$ with $\|Y\|=1$ 
satisfying the following conditions:
\begin{itemize}
\item[(i)]  $\nabla^L_Z Y=0$,
\item[(ii)]  $\alpha_g(Y,Z)=0$,
\item[(iii)] $\<(\gamma I+\textup{Hess}\,\gamma)Y,Z\>_L=0$.
\end{itemize}
Then the orthogonal complement of $\textup{span}\{jY\}$ is a totally geodesic 
distribution on the open subset $V$ of regular points of $\psi$
endowed with the induced metric, where $j=j(y,w)\colon T_yL\to \Delta^\perp(y,w)$ 
is the isometry given by part $(iii)$ of Proposition~\ref{gaussprop} between $T_yL$ 
and the orthogonal complement of $\Delta(y,w)=N_gL(y)$ in $T_{(y,w)}\Lambda$.
Moreover, $\psi$ is ruled if $Y$ and $Z$ are linearly dependent and a 
partial tube over a curve with fiber a flat hypersurface if $Y$ and $Z$
are orthogonal.

Conversely, if $V$ carries a  totally geodesic distribution of codimension 
one then there exist $Y,Z\in\mathfrak{X}(L)$ with $\|Y\|=1$ 
satisfying conditions $(i)$ to $(iii)$. 
\end{lemma}

\proof Given $(y,w)\in V$, let $P_w(y)$ be the invertible endomorphism 
of $T_yL$ given by (\ref{eq:pw}). Since conditions  $(ii)$ and $(iii)$ are 
satisfied, then
\be\label{pw}
\<P_wZ,Y\>_L=0
\ee
for all $w\in\Gamma(N_gL)$. Thus $P_wZ$ is collinear with $X$ for all  
$w\in \Gamma(N_gL)$, where $X\in\mathfrak{X}(L)$ is a unit vector field 
orthogonal to $Y$. Hence, for any $(y,w)\in V$ the vector
 $P^{-1}_wX(y)$ is collinear with $Z(y)$.    By $(ii)$ we have
$$
\<A_{w'}Z,Y\>_L=\<\alpha_g(Y,Z),w'\>=0
$$
for all $w'\in\Gamma(N_gL)$. Using  \eqref{cxiaxi} and the fact that
$j$ is an isometry we obtain
\be\label{eq:awpw2}
\<\nabla^V_{jX}jY,w'\>_V=\<C_{w'}jX,jY\>_V=\<A_{w'}P_w^{-1}X,Y\>_L=0
\ee
at any $(y,w)\in V$ and
for all $w'\in N_gL(y)$. On the other hand, using \eqref{eq:conngpsi} 
we obtain
\be\label{eq:awpw3}
\<\nabla^V_{jX}jY,jX\>_V=\<\nabla^L_{P_w^{-1}X}Y,X\>_L=0
\ee
where in the last equality we have used condition $(i)$
and the fact that $P^{-1}_wX(y)$ is collinear with $Z(y)$.
We conclude from \eqref{eq:nxijx}, \eqref{eq:awpw2} 
and \eqref{eq:awpw3} that the distribution orthogonal to  $jY$  is totally 
geodesic.

We prove next the last assertion of the direct statement. Suppose first that 
$Y$ and $Z$ are linearly dependent. Without loss of generality, we may assume that 
$Z=Y$. As before, let $X\in\mathfrak{X}(L)$ be a unit vector field orthogonal 
to $Y$. By \eqref{eq:sffpsi}, at any $(y,w)\in V$  we have
$$
\<AjX,jX\>_V=-\<jP_w^{-1}X,jX\>_V=-\<P_w^{-1}X,X\>_L.
$$
On the other hand, since \eqref{pw} holds, then
$$ 
\<P_wY,Y\>_L=0,
$$
hence 
$$ 
P_wY=\<P_wY,X\>_LX.
$$
Therefore 
$$
Y=\<P_wY,X\>_LP_w^{-1}X,
$$
which implies using (\ref{eq:sffpsi}) that 
$$
0=\<P_w^{-1}X,X\>_L=-\<AjX,jX\>_V,
$$
and we conclude that $\psi$ is ruled.  

Assume now that $Z$ and $Y$ are everywhere orthogonal. Arguing in a similar way, 
we see that $AjZ$ is collinear with  $jZ$, which implies from Proposition \ref{pt2}
that $\psi$ is a partial tube over a curve $\gamma$, with fiber a hypersurface 
$N^{n-1}$ of the affine normal space to $\gamma$  at some point. Since $\psi$ has 
rank two, $N^{n-1}$ must be flat.  

Next we prove the converse. If  $\tilde Y\in\Gamma(\Delta^\perp)$ is a 
unit vector field such that the distribution $\{\tilde Y\}^\perp$ is totally 
geodesic, then
$$
\nabla^V_{\tilde X}\tilde Y=0
$$
for any $\tilde X\in\mathfrak{X}(V)$ such that $\<\tilde X,\tilde Y\>_V=0$.
In particular,
$$
\nabla^V_T\tilde Y=0
$$
for all $T\in\Gamma(\Delta)$. Hence there exists $Y\in\mathfrak{X}(L)$ such
that $\tilde Y=j Y$. From
$$
\nabla^V_{jX}jY=0
$$
for $X\in\mathfrak{X}(L)$ such that $\<X,Y\>_L=0$, using
\eqref{cxiaxi} and \eqref{eq:conngpsi} we obtain, respectively, 
\be\label{eq:awpw2x}
0=\<\nabla^V_{jX}jY, w'\>_V=\<C_{w'}jX, jY\>_V=\<A_{w'}P_w^{-1}X,Y\>_L
\ee
and
\be\label{eq:w1}
\nabla^L_{P_w^{-1}X}Y=0
\ee
for any $(y,w)\in V$ and for all $w'\in N_gL(y)$. 

Suppose that at $y\in L^2$ we have 
$$
T_yL=\textup{span}\{P_w^{-1}X(y) : w\in N_gL(y)\;\mbox{with}\;(y,w)\in V\}.
$$
Then this also holds in a neighborhood $W$ of $y$, and then \eqref{eq:w1}
implies that $Y$ is parallel on $W$. We conclude that $W$ is flat.
On the other hand, by \eqref{eq:awpw2x} we have 
$$
\<P_w^{-1}X(y),A_{w'}Y\>_L=0
$$
for any $w\in N_gL(y)$ with $(y,w)\in V$, and for all $w'\in N_gL(y)$.  
It follows that $A_{w'}Y=0$ for all
$w'\in \Gamma(N_gW)$. But then $Y$  belongs to the relative nullity
space of $g$ at any point of $W$, and hence the Gauss curvature of $L^2$
is equal to one at any point of $W$, a contradiction.

It follows that there exists  $Z\in\mathfrak{X}(L)$  such that,
for all $y\in L^2$ and for all $w\in N_gL(y)$ such that $P_w$ is invertible,
the vector $P_w^{-1}X(y)$ is collinear with $Z(y)$.  This implies that
$$
\<P_wZ,Y\>_L=0.
$$
On the other hand,  equation \eqref{eq:awpw2x} gives
$$
\<A_wZ,Y\>_L=0
$$
for all $y\in L^2$ and $w\in N_gL(y)$, hence
$$
\<(\gamma I+\textup{Hess}\,\gamma)Y,Z\>_L=0.
$$
Finally, from \eqref{eq:w1} it follows that $\nabla^L_ZY=0$.
\vspace{1ex}\qed

Recall that a  smooth  distribution $D$ 
on a Riemannian manifold $M^n$ is said to be \emph{curvature invariant} if
$$
R(X,Y)Z \in D\;\;\mbox{for all}\;\; X,Y,Z\in D
$$
where $R$ denotes the curvature tensor of $M^n$.  
The next result was obtained in \cite{drt}. 

\begin{proposition}\po\label{prop:curv}
Let $f\colon M^n\to\Q_\e^{n+1}$ be an oriented hypersurface carrying a 
curvature invariant distribution  $D$ of rank $k$.
Then one of the following possibilities holds pointwise:
\begin{itemize}
\item[(i)] $A(D)\subset D^\perp$,
\item[(ii)] $A(D)\subset D$,
\item[(iii)] $\textup{rank}\,D\cap \Delta=k-1$.
\end{itemize}
\end{proposition} 

\noindent\emph{Proof of Theorem \ref{thm:main}:} 
By Proposition \ref{prop:curv},  the shape operator $A$ of $f$ 
satisfies pointwise either of conditions $(i)$ to $(iii)$ in that 
result. In any oriented open subset where $A$ satisfies condition 
$(ii)$, $f$ is locally a partial tube over a curve by 
Proposition \ref{pt2}. On the other hand, $f$ is ruled on an oriented 
open subset where condition $(i)$ is satisfied. Finally, on an oriented 
open subset $U$ where $A$ satisfies condition $(iii)$ 
and neither of conditions $(i)$ or $(ii)$,  $f$ has rank two since 
$M^n$ is nowhere flat.
By Lemma \ref{totgeofol}, $f$ can be parametrized in terms of 
the Gauss parametrization by a pair $(g,\gamma)$, with 
$g\colon L^2\to\Sf^n$, $n\geq 3$, and $\gamma\in C^{\infty}(L)$, 
for which there exist $Y,Z\in\mathfrak{X}(L)$ with $\|Y\|=1$ 
satisfying conditions $(i)$ to $(iii)$ in that lemma. 

Since $f$ is nowhere ruled or a partial tube over a curve on $U$, 
the vector fields $Y, Z$ are linearly independent and 
nowhere orthogonal on $U$. Therefore one of the following possibilities 
holds locally on $U$: $g$ is totally geodesic and $\gamma\in C^{\infty}(L)$
is any smooth function such that $P=\gamma I+\textup{Hess}\,\gamma$ is 
invertible, and hence $f$ is cylindrical surfacelike, or $g(L^2)$ is 
contained in an umbilical $\Sf^3\subset \Sf^n$ and $\gamma=0$, in 
which case $f$ is conical surfacelike, or the pair $\{Y,Z\}$ satisfies
conditions (\ref{eq:typeDcond}) for $g$ and $\gamma$ satisfies the
same differential equation (\ref{eq:pos}) as the position vector of $g$, 
by virtue of condition $(iii)$, and hence $(g,\gamma)$ is a pair of type 
$D$. The converse statement follows from the
direct statement of Lemma \ref{totgeofol}. \qed

\begin{remark}\po {\em Theorem \ref{thm:main} can be used to 
derive the result in \cite{drt} referred to in the introduction
for foliations with complete leaves, once we first prove, as in 
\cite{drt}, that if $U\subset M^n$ is an oriented open subset
at each point of which  neither of conditions $(i)$ or $(ii)$ 
in Proposition \ref{prop:curv} is satisfied, then any unit speed 
geodesic starting at a point of $U$ and tangent to $\Delta$ at 
that point remains in $U$ for any value of the parameter. 

  In fact, after proving  this, since $f$ has rank two on $U$ and the leaves 
of $\Delta$ are complete therein, one can globally parametrize
$f$ on $U$ by means of the Gauss parametrization $\psi$ given 
by (\ref{gpara}) in terms of a pair $(g, \gamma)$ of type $D$, 
and then we may identify $M^n$ with the subset $V$ of regular points 
of $\psi$, endowed with the induced metric, and $f$ with $\psi|_V$.
Then, given $(y,w)\in U\subset V$, any geodesic $\gamma$ starting 
at $(y,w)$ and tangent to $\Delta(y,w)=N_gL(y)$ at $(y,w)$ is given 
by $t\mapsto (y,tw')$ for some $w'\in N_gL(y)$. 

Let $Y,Z\in \mathfrak{X}(L)$ satisfy conditions $(i)$ to $(iii)$ 
in the definition of surfaces of type $D$. Since $(g,\gamma)$ is 
a pair of type $D$, the matrix of the endomorphism 
$$
P_{tw'}=\gamma(y)I+\textup{Hess}\,\gamma(y)-tA_{w'}
$$ 
of $T_yL$ with respect to the basis $\{Y(y),Z(y)\}$ is a diagonal 
matrix whose diagonal elements are linear in $t$, and hence have zeros 
unless the shape operator $A_{w'}$ is identically zero. 

  In summary, there always exist points $(y,tw')$ along $\gamma$ at which $\psi$ fails 
to be regular, unless $A_{w'}$ is identically zero. Thus the latter 
possibility must hold for all $(y,w)\in U$ and for all $w'\in N_gL(y)$. 
But this implies that $\psi(U)$ is a surfacelike strip, which is 
ruled out by assumption. The proof then proceeds as in \cite{drt} by 
showing that this implies $\psi$ to be either
ruled or a partial tube over a curve.}
\end{remark}

\section[The case of nonflat space forms]
{The case of nonflat space forms}

 As in the Euclidean case, among the hypersurfaces 
$f\colon M^n\to\Q_\e^{n+1}$ in nonflat ambient space forms, that is, 
$\e\in\{1,-1\}$, that carry totally geodesic distributions of 
codimension one, one has the ruled hypersurfaces, that is, the 
hypersurfaces carrying a smooth foliation by (open subsets of) 
$(n-1)$-dimensional totally geodesic submanifolds of $\Q_\e^{n+1}$, 
and the  partial tubes over curves described in Section $1.1$. 
There are also two other families of trivial examples.
\vspace{1ex} 
 
\noindent $(1)$ The first family  is that of isometric immersions 
$f\colon U\to\Q_e^{n+1}$ of open subsets $U\subset\Q_\e^n$, for one 
may consider any foliation of $U$  by (open subsets of) totally 
geodesic hypersurfaces of constant sectional curvature $\e$.  
Notice that these examples may have complete leaves if $\e=-1$.
\vspace{1ex} 

\noindent $(2)$ The second family of examples consists of generalized
cones. Let $g\colon L^2\to\Q_c^3$, $c\geq \e$, be an isometric immersion, 
and let  $i\colon\Q_c^3\to\Q_\e^{n+1}$ be an umbilical inclusion. 
Thus, the normal bundle of $\tilde g=i\circ g$ splits as 
$$
N_{\tilde g}L=i_*N_gL\oplus N_i\Q_c^3.
$$ 
We regard $\Lambda=N_i\Q_c^3$ as a subbundle of 
$N_{\tilde g}L$ and define $f\colon\Lambda\to\Q_\e^{n+1}$ by 
$$
f(x,v)=\exp_{g(x)}v,
$$
where $\exp$ is the exponential map of $\Q_\e^{n+1}$. Then $f$ is called 
the \emph{generalized cone over}~$g$. Here again, if $\e=-1$ and $c\leq 0$, 
then generalized cones yield examples of hypersurfaces carrying totally 
geodesic foliations of codimension one with complete leaves.
\vspace{1ex}
 
The local classification of  hypersurfaces of nonflat space forms that 
carry a totally geodesic  foliation of codimension one is as follows.

\begin{theorem}\po\label{thm:main2}
Let $f\colon M^n\to\Q_\e^{n+1}$, $n\geq 3$ and $\e\in \{1,-1\}$, 
be an isometric immersion of a Riemannian manifold without points 
where all sectional curvatures are  equal to~$\e$. If $M^n$ carries 
a totally geodesic  foliation of codimension one, then there exists an 
open dense subset of $M^n$ where  $f$ is locally either ruled, a partial 
tube over a curve, a generalized cone or given in terms of the Gauss 
parametrization by a surface $g\colon L^2\to\Sf_\mu^{n+1}$ of type~$D$, 
\mbox{where $\mu=(1-\e)/2$}. 
\end{theorem}

Similarly to the case of Euclidean hypersurfaces, the proof relies on the 
following version of Lemma \ref{totgeofol}.

\begin{lemma}\po\label{totgeofol4}
Let $g\colon L^{2}\to\Sf_\mu^{n+1}$ be a surface and let
$\Lambda_\e$ and  $\psi\colon\Lambda_\e\to\Q_\e^{n+1}$ be 
given by (\ref{gammae}) and (\ref{gpara2}), respectively.  
Assume that there exist $Y,Z\in\mathfrak{X}(L)$ with $\|Y\|=1$
satisfying conditions $(i)$ and $(ii)$ in Lemma \ref{totgeofol}.  
Then the orthogonal complement of $\textup{span}\{jY\}$ is a 
totally geodesic distribution on the open subset $V$ of regular 
points of $\psi$, endowed with the induced metric.
Here  $j=j(y,w)\colon T_yL\to \Delta^\perp(y,w)$ is the isometry 
given by part $(iv)$ of Proposition \ref{gauss3}  between $T_yL$ 
and the orthogonal complement of
$$
\Delta(y,w)=\{w'\in N_gL(y):\<w', w\>=0\}\subset T_{(y,w)}\Lambda_\e.
$$ 

Conversely, if $V$ carries a  totally geodesic distribution of codimension 
one then there exist $Y,Z\in\mathfrak{X}(L)$ with $\|Y\|=1$
satisfying conditions $(i)$ and $(ii)$ in Lemma \ref{totgeofol}. 
\end{lemma}

\proof  It follows from $(ii)$  that 
$$
\<A_wZ,Y\>_L=0
$$
for all $w\in \Gamma(N_gL)$.  Therefore, if $X\in\mathfrak{X}(L)$ 
is a unit vector field orthogonal to $Y$, then $A_wZ$ is collinear 
with $X$ for all $w\in \Gamma(N_gL)$. Hence, whenever $A_w$ is invertible,
we see that $A^{-1}_wX$ is collinear 
with $Z$.  
It follows from \eqref{cxiaxi2} that 
\be\label{eq:awpw2c}
\<\nabla^V_{jX}jY,w'\>_V=\<C_{w'}jX,jY\>_V=\<A_{w'}A_w^{-1}X,Y\>_L=0
\ee
for any $(y,w)\in V$ and for all $w'\in\Delta(y,w)$. On the other 
hand, by \eqref{eq:conngpsi2} we have
\be\label{eq:awpw3b}
\<\nabla^V_{jX}jY,jX\>_V=\<\nabla^L_{A_w^{-1}X}Y,X\>_L=0
\ee
by $(i)$ and the fact that $A^{-1}_wX$ is collinear with $Z$.
 We conclude from \eqref{eq:nxijx}, \eqref{eq:awpw2c} 
and \eqref{eq:awpw3b} that the distribution orthogonal to $jY$ is 
totally geodesic.

We now prove the converse. If  $\tilde Y\in \Gamma(\Delta^\perp)$ is a 
unit vector field such that the distribution $\{\tilde Y\}^\perp$ is totally 
geodesic, then
$$
\nabla^V_{\tilde X} \tilde Y=0
$$
for any $\tilde X\in\mathfrak{X}(V)$ such that $\<\tilde X,\tilde Y\>_V=0$.
In particular,
$$
\nabla^V_T \tilde Y=0
$$
for all $T\in \Gamma(\Delta)$, hence there exists 
$Y\in\mathfrak{X}(L)$ such that $\tilde Y=j Y$. From
$$
\nabla^V_{jX}jY=0
$$
for $X\in\mathfrak{X}(L)$ such that $\<X,Y\>_L=0$,  using
\eqref{cxiaxi2} and \eqref{eq:conngpsi2} we obtain, respectively, 
\be\label{eq:awpw2j}
0=\<\nabla^V_{jX}jY, w'\>_V=\<C_{w'}jX, jY\>_V=\<A_{w'}A_w^{-1}X,Y\>_L
\ee
and
\be\label{eq:w1b}
\nabla^L_{A_w^{-1}X}Y=0
\ee
for any $(y,w)\in V$ and for all $w'\in\Delta(y,w)$.

Suppose that at $y\in L^2$ we have
$$
T_yL=\textup{span}\{A_w^{-1}X(y): w\in N_gL(y)\;\mbox{with}\;(y,w)\in V\}.
$$
Then this also holds in a neighborhood $W$ of $y$, and then  $Y$ 
is parallel on $W$ by \eqref{eq:w1b}, which implies that $V$ is flat.
On the other hand, by \eqref{eq:awpw2j} we have
$$
\<A_w^{-1}X,A_{w'}Y\>_L=0
$$
for all $w, w'\in N_gL(y)$ such that $A_w$ is invertible and $\<w', w\>=0$.
 It follows that $A_{w'}Y=0$ for all
$w'\in \Gamma(N_gW)$. But then $Y$ belongs to the relative nullity
space of $g$ at any point of $W$, and hence the Gauss curvature of $L^2$
is equal to one at any point of $W$, a contradiction.

It follows that there exists  $Z\in\mathfrak{X}(L)$  such that,
for all $y\in L^2$ and for all $w\in N_gL(y)$ such that $A_w$ is 
invertible, the vector $A_w^{-1}X(y)$ is collinear with $Z(y)$. 
Eq. \eqref{eq:awpw2c} and the fact that $X$ and $Y$ are orthogonal 
then imply that
$$
\<A_wZ,Y\>_L=0
$$
for all $w\in N_gL(y)$. Finally, from \eqref{eq:w1b} we obtain  
$\nabla^L_ZY=0$.
\vspace{1ex}\qed

Our last result shows that  the only hypersurfaces of  $\Sf^{n+1}$,  
$n\geq 3$, that may carry a totally geodesic distribution of codimension 
one with complete leaves are  partial tubes over curves. 

\begin{theorem}\po\label{thm:main3}
Let $f\colon  M^n\to\Sf^{n+1}$, $n\geq 3$, be an isometric immersion of 
a Riemannian manifold without  points where all sectional curvatures are 
equal to $1$. If $M^n$ carries a totally geodesic foliation of codimension 
one with complete leaves, then $\tilde{f}=f\circ\pi$ is a partial tube over a 
curve, where $\pi:\tilde{M}\to M^n$ is the universal cover space of $M^n$.
\end{theorem}

\proof  First we prove that there cannot exist any open subset 
$U\subset M^n$ where neither of conditions $(i)$ or $(ii)$ in 
Proposition \ref{prop:curv} hold.
Suppose otherwise and let $U\subset M^n$ be the maximal open subset where 
neither of conditions $(i)$ or $(ii)$. Then condition $(iii)$ in Proposition 
\ref{prop:curv} holds and  $\Delta\subset D$ at any point of $U$. 
Choose unit vector fields  $Y\in\Gamma(D^\perp)$ 
and $X\in\Gamma(D\cap\Delta^\perp)$ along $U$. 
Then there exist $\beta,\mu,\rho\in C^{\infty}(U)$, with $\rho\mu\neq 0$ 
at any point of $U$, such that
$$
AY=\beta Y+\mu X\;\;\;\mbox{and}\;\;\;AX=\mu Y+\rho X.
$$
Since $D$ and $\Delta$ are both totally geodesic and 
$\Delta\subset D$, we have
$$
\nabla_XY=\nabla_TY=\nabla_TX=0
$$
for all $T\in \Gamma(\Delta)$. The Codazzi equation for
$\left(X,T\right)$ yields
\be\label{eq:cod1}
T(\rho)=\rho\,\<\nabla_XX,T\>\;\;\;\mbox{and}\;\;\;
T(\mu)=\mu\,\<\nabla_XX,T\>
\ee
for all $T\in\Gamma(\Delta)$, whereas the  Gauss equation 
for $\left(X,T,S,X\right)$ gives
\be\label{eq:gauss1}
T\<\nabla_XX,S\>=\<T,S\>+\<\nabla_XX,\nabla_T S\>
+\<\nabla_XX,T\>\<\nabla_XX,S\>
\ee
for all $T,S\in\Gamma(\Delta)$. 

 From (\ref{eq:cod1}) we see that $\mu=\varphi\rho$ with $T(\varphi)=0$ 
for all $T\in\Gamma(\Delta)$. 
It follows that the leaves of $\Delta$ through points of $U$ are
complete and that condition $(iii)$ remains valid along them. Otherwise, 
there would exist a  geodesic tangent to some $T\in \Gamma(\Delta)$ of 
unit length that would reach a point of the boundary of $U$. Since such a 
point is the limit of a sequence of points where either of conditions $(i)$ 
or $(ii)$  holds, either  $\mu$ or $\rho$  
vanishes at that point, and hence both $\mu$ and $\rho$ by the above.
But this possibility is ruled out by the assumption that there are no 
points in $M^n$  where the sectional curvatures are all equal to $1$.

Now, for any complete geodesic with unit tangent vector $T$ in a leaf 
of  $\Delta$, by (\ref{eq:gauss1}) the function   
$\lambda=\<\nabla_XX,T\>$ satisfies the differential equation
\be\label{eq:lambda}
T(\lambda)=1+\lambda^2,
\ee
and thus $\arctan\lambda(t)=t$, a contradiction with the fact that 
$\lambda$ is defined for any value of the parameter. 

It follows that $U$ is empty, and hence either of conditions $(i)$ 
or $(ii)$ must hold at any point of $M^n$. Let $S_1$ (respectively, $S_2$) 
be the subset of $M^n$ where condition $(i)$ (respectively, condition $(ii)$) 
is satisfied.  Since both $S_1$ and $S_2$ are closed and 
$M^n=S_1\cup S_2$, any point $x\in\partial S_1$ belongs to $S_1\cap S_2$, 
hence all sectional curvatures are equal to $\e$ at $x$.
It follows from our assumption that either $M^n=S_1$ or $M^n=S_2$. 

In the first case, $f$ is a ruled hypersurface. But then $f$ has still 
rank two, $\Delta\subset D$ and (\ref{eq:gauss1}) 
still holds, and hence also (\ref{eq:lambda}) along  each complete geodesic 
with unit tangent vector $T$ in a leaf of $\Delta$, which leads to a 
contradiction as before.

We conclude that $M^n=S_2$. Since the leaves of $D$ are complete, then
the universal covering space of $M^n$ is isometric to a twisted product 
$\tilde{M}^n=\tilde{M}_0^{n-1}\times_\rho\mathbb{R}$, where 
$\rho\in C^\infty(\tilde{M})$ (see Theorem $2.7$ of \cite{to}), with 
the leaves of the lifting $\tilde{D}$ of $D$ corresponding to the 
slices $\tilde{M}_0^{n-1}\times\{t\},\ t\in\mathbb{R}$. Moreover, 
the fact that $M^n=S_2$ means that the vector field $\d/\d t$ 
tangent to the second factor is a principal direction of $\tilde f$ 
at any point. It follows from Proposition~\ref{pt2} that $\tilde f$ 
is a partial tube over a curve.\vspace{1ex}\qed

\noindent {\bf Acknowledgment.\/}  We are very grateful to the
anonymous referee for his careful reading of the paper and for several 
useful comments and suggestions.

{\renewcommand{\baselinestretch}{1}
\hspace*{-30ex}\begin{tabbing}
\indent \= IMPA  \hspace{30ex} Universidade de S\~ao Paulo \\
\>  Estrada Dona Castorina, 110 \hspace{7ex}
Av. Trabalhador S\~ao-Carlense 400 \\
\> 22460-320 --- Rio de Janeiro
\hspace{7.5ex} 13560-970 --- S\~ao Carlos  \\
\> Brazil\hspace{31ex} Brazil\\
\> marcos@impa.br  \hspace{20ex}
tojeiro@icmc.usp.br
\end{tabbing}}
\end{document}